\documentclass[12pt,letterpaper]{amsart}
\usepackage{epic,eepic,latexsym, amssymb, amscd, amsfonts, xypic}

 \newlength{\baseunit}               
 \newcount{\numlines}                
\newcommand{\Z}{\mathbb{Z}}

\newcommand{\C}{\mathbb{C}}

\newcommand{\proj}{\mathbb P}
\newcommand{\oh}{{\mathcal{O}}}

\newcommand{\cm}{{\mathcal{M}}}
\newcommand{\fM}{{\mathfrak{M}}}
\newcommand{\cmbar}{\overline{\cm}}

\newcommand{\be}{\beta}

\newcommand{\De}{\Delta}

\newcommand{\Bl}{\operatorname{Bl}}

\newcommand{\secretnote}[1]{}

\newcommand{\lremind}[1]{{}}

\begin{document}
\pagestyle{plain}

\title{\large{
A natural smooth compactification of the space of
elliptic curves in projective space
}}
\author{Ravi Vakil and Aleksey Zinger}
\date{July 14, 2006.}
\begin{abstract}
  The space of smooth genus $0$ curves in projective space has a
  natural smooth compactification: the moduli space of stable maps,
  which may be seen as the generalization of the classical space of
  complete conics.  In arbitrary genus, no such
natural  smooth model is expected, as the space satisfies ``Murphy's Law''.
  In genus $1$, however, the situation remains beautiful.  We give a
  natural smooth compactification of the space of elliptic curves in
  projective space, and describe some of its properties.  This space
  is a blow up of the space of stable maps.  It can be interpreted as
  blowing up the most singular locus first, then the next most
  singular, and so on, but with a twist --- these loci are often
  entire components of the moduli space.  We give a number of
  applications in enumerative geometry and Gromov-Witten theory.
The proof that this construction indeed gives a desingularization
will appear in \cite{vz}.
\end{abstract}

\maketitle

{\parskip=12pt 


The moduli space of stable maps $\cmbar_{g,k}(X, \be)$ to a complex
projective manifold $X$ (where $g$ is the genus, $k$ is the number of
marked points, and $\be \in H_2(X, \Z)$ is the image homology class)
is the central tool and object of study in Gromov-Witten theory.  We
consider this space as a Deligne-Mumford stack.  The open subset
corresponding to maps from smooth curves is denoted $\cm_{g,k}(X,
\be)$.

The protean example is $\cmbar_{0,k}(\proj^n,d)$. 
This space is wonderful in essentially all ways:
it is irreducible, smooth, contains $\cm_{0,k}(\proj^n,d)$
as a dense open subset.  The boundary
$$\De := \cmbar_{0,k}(\proj^n,d) \setminus
\cm_{0,k}(\proj^n,d)$$ is normal crossings.  The divisor theory is
fully understood, and combinatorially tractable,
 \cite{p}.
In some sense, this should be seen as the natural generalization
of the space of complete conics
compactifying the space of smooth conics.

It is natural to wonder if such a beautiful structure
exists in higher genus.  In arbitrary genus,
however, there is no reasonable hope:
even the interior $\cm_g(\proj^n,d)$ is
badly behaved in general.
More precisely, $\cm_g(\proj^n,d)$ (as $g$, $n$, and $d$ vary)
is arbitrary singular in 
a well-defined sense --- it can have essentially any
singularity, and can have components of various dimension meeting in
various ways with various nonreduced structures \cite{v}.
There is no reasonable hope of describing a desingularization,
as this would involve describing a resolution of
singularities.

In genus one, however, the situation remains remarkably beautiful.
Although $\cmbar_{1,k}(\proj^n,d)$ in general has many components, it
is straightforward to show that $\cm_{1,k}(\proj^n,d)$ is irreducible
and smooth.  Let $\cmbar_{1,k}^0(\proj^n,d)$ be the closure of this
open subset (the ``main component'' of the moduli space).

In the paper \cite{vz}, we will show that there is a natural desingularization
of this main component $$\widetilde{\cm}_{1,k}(\proj^n,d)
  \rightarrow \cmbar_{1,k}^0(\proj^n,d).$$
This desingularization has several desirable properties. 
\begin{itemize}
\item It leaves the interior $\cm_{1,k}(\proj^n,d)$ unchanged.
\item The boundary $\widetilde{\cm}_{1,k}(\proj^n,d) \setminus \cm_{1,k}(\proj^n,d)$ 
is simple normal crossings, with an explicitly described
normal bundle.  
\item The points of the boundary have an explicit geometric interpretation.
\item The desingularization can be interpreted as blowing up
``the most singular locus'', then ``the next most singular locus'',
and so on, but with an unusual twist.
\item The divisor theory is explicitly describable, and the intersection
theory is tractable.  (For example, one can compute the top
intersection of any combination of divisors using \cite{z}.)
\item The compactification is natural in the following senses.
\begin{enumerate}
\item[(i)] The desigularization is equivariant --- it behaves well
  with respect to the symmetries of $\proj^n$.  Hence we can apply
  Atiyah-Bott localization to this space --- not just in theory, but
  in practice.
\item[(ii)] It behaves well with respect to the
inclusion $\proj^m \hookrightarrow \proj^n$.
\item[(iii)] It behaves well with respect to the marked points
(forgetful maps; $\psi$-classes; etc.).
\item[(iv)] 
Consider the universal map $$\xymatrix{ \mathcal{C} \ar[r]^{\pi} \ar[d]^{\rho} &
\proj^n \\ \cmbar_{g,k}(\proj^n,d)}.$$
An important sheaf in Gromov-Witten theory is $\rho_* \pi^* \oh_{\proj^n}(a)$.
When $g>0$, this is not a vector bundle, which causes
difficulty in computation and theory.  However, in genus $1$, 
``resolving $\cmbar_{1,k}^0(\proj^n,d)$ also resolves this sheaf'':
when the sheaf is pulled back to the desingularization, 
it ``becomes'' a vector bundle.  More precisely:  it contains
a natural vector bundle, and is isomorphic to it on the interior.
This  vector bundle is explicitly describable.
\end{enumerate}
\end{itemize}

We think it is interesting that such a natural naive approach as we
describe below actually works, and yields a desingularization with
these nice properties.  For example, if $n>2$, this desingularization
can be interpreted as a natural compactification of the Hilbert scheme
of smooth degree $d$ curves in projective space, and thus could be
seen as the genus $1$ version of the complete conics.

This construction also has a number of applications: 
\begin{itemize}
\item enumerative geometry of genus $1$ curves via localization
(extending results of \cite{crelle}, for example adding tangencies).
\item Gromov-Witten invariants in terms of enumerative invariants  \cite{reduced}.
\item the Lefschetz hyperplane property:  effective computation of Gromov-Witten invariants  of complete intersections \cite{lz2} (see also \cite{lz1} 
for the special case of the quintic threefold).
\item an algebraic version of ``reduced'' Gromov-Witten invariants
in symplectic geometry \cite{reduced}. 
\item  an approach to the physicists' prediction \cite{bcov} of
genus $1$ Gromov-Witten invariants \cite{zip}.
\end{itemize}

Before giving the construction, we motivate it by describing the
geography of $\cmbar_{1,k}(\proj^n,d)$ It is straightforward to show
that $\cmbar_{1,k}(\proj^n, d)$ is nonsingular on the locus where
there is no contracted genus $1$ (possibly nodal) curve (for example,
the proof of \cite[Prop.~4.21]{crelle} applies).

\noindent {\bf Example:  plane cubics.}

We first consider the case of $\cmbar_1(\proj^2,3)$, see
Figure~\ref{cubics}.  The main component generically corresponds to
smooth plane cubics, which has dimension $9$.  This is depicted in the
upper-central panel of the figure.  The remaining components must all
contain a contracted genus $1$ curve, and we enumerate the
possibilities.

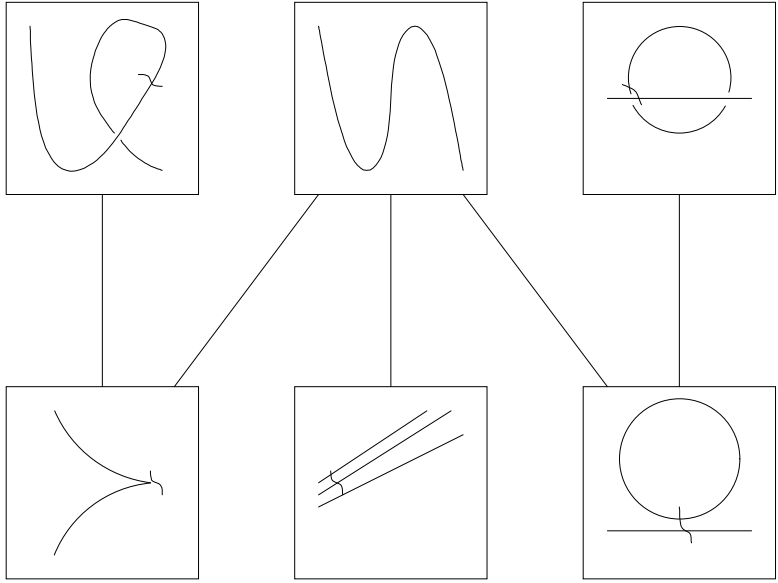
\begin{figure}
\setlength{\unitlength}{0.00083333in}
\begingroup\makeatletter\ifx\SetFigFont\undefined%
\gdef\SetFigFont#1#2#3#4#5{%
  \reset@font\fontsize{#1}{#2pt}%
  \fontfamily{#3}\fontseries{#4}\fontshape{#5}%
  \selectfont}%
\fi\endgroup%
{\renewcommand{\dashlinestretch}{30}
\begin{picture}(4824,3639)(0,-10)
\path(3855,3099)(3856,3099)(3863,3096)
	(3877,3091)(3892,3084)(3907,3078)
	(3918,3072)(3927,3065)(3935,3058)
	(3941,3051)(3947,3041)(3953,3028)
	(3960,3013)(3966,2996)(3971,2981)
	(3974,2973)(3974,2972)
\put(981.643,-105.857){\arc{1442.455}{3.5221}{4.6157}}
\put(4212.783,3141.900){\arc{640.201}{2.8209}{6.5679}}
\put(4210.085,3123.805){\arc{657.614}{0.5040}{2.6514}}
\put(4212,762){\ellipse{750}{750}}
\path(12,3612)(1212,3612)(1212,2412)
	(12,2412)(12,3612)
\path(1812,3612)(3012,3612)(3012,2412)
	(1812,2412)(1812,3612)
\path(3612,3612)(4812,3612)(4812,2412)
	(3612,2412)(3612,3612)
\path(3612,1212)(4812,1212)(4812,12)
	(3612,12)(3612,1212)
\path(1812,1212)(3012,1212)(3012,12)
	(1812,12)(1812,1212)
\path(12,1212)(1212,1212)(1212,12)
	(12,12)(12,1212)
\path(1962,612)(2637,1062)
\path(1962,537)(2787,1062)
\path(1962,462)(2862,912)
\path(3762,312)(4662,312)
\path(3762,3012)(4662,3012)
\path(1962,2412)(1062,1212)
\path(2412,2412)(2412,1212)
\path(2862,2412)(3762,1212)
\path(4212,2412)(4212,1212)
\path(612,2412)(612,1212)
\path(1962,3462)(1962,3460)(1963,3456)
	(1964,3448)(1967,3435)(1970,3418)
	(1974,3395)(1979,3367)(1985,3335)
	(1991,3297)(1999,3256)(2008,3212)
	(2017,3165)(2027,3117)(2037,3068)
	(2048,3019)(2059,2970)(2071,2922)
	(2083,2876)(2096,2831)(2109,2788)
	(2123,2748)(2137,2710)(2153,2675)
	(2169,2644)(2186,2617)(2204,2594)
	(2223,2577)(2242,2566)(2262,2562)
	(2282,2566)(2301,2577)(2319,2594)
	(2335,2616)(2350,2643)(2363,2673)
	(2373,2706)(2382,2741)(2390,2777)
	(2396,2815)(2400,2853)(2404,2893)
	(2407,2932)(2410,2972)(2412,3012)
	(2414,3052)(2417,3092)(2420,3131)
	(2424,3171)(2428,3209)(2434,3247)
	(2442,3283)(2451,3318)(2461,3351)
	(2474,3381)(2489,3408)(2505,3430)
	(2523,3447)(2542,3458)(2562,3462)
	(2582,3458)(2601,3447)(2620,3430)
	(2638,3407)(2655,3380)(2671,3349)
	(2687,3314)(2701,3276)(2715,3236)
	(2728,3193)(2741,3148)(2753,3102)
	(2765,3054)(2776,3005)(2787,2956)
	(2797,2907)(2807,2859)(2816,2812)
	(2825,2768)(2833,2727)(2839,2689)
	(2845,2657)(2850,2629)(2854,2606)
	(2857,2589)(2860,2576)(2861,2568)
	(2862,2564)(2862,2562)
\path(162,3462)(162,3460)(162,3456)
	(162,3449)(163,3437)(163,3420)
	(164,3399)(165,3373)(166,3342)
	(168,3307)(170,3268)(173,3226)
	(176,3182)(179,3136)(183,3089)
	(188,3042)(193,2995)(199,2949)
	(206,2904)(215,2860)(224,2817)
	(234,2777)(247,2738)(260,2702)
	(276,2669)(294,2638)(314,2612)
	(336,2590)(360,2573)(387,2562)
	(415,2558)(444,2560)(473,2567)
	(501,2578)(528,2592)(554,2609)
	(578,2627)(600,2647)(622,2668)
	(641,2689)(660,2711)(678,2734)
	(695,2756)(711,2779)(727,2802)
	(743,2826)(758,2850)(774,2875)
	(791,2900)(808,2926)(825,2953)
	(843,2981)(862,3011)(881,3041)
	(901,3073)(921,3105)(940,3138)
	(958,3172)(974,3205)(987,3237)
	(997,3269)(1004,3299)(1007,3326)
	(1007,3350)(1004,3370)(1000,3387)
	(993,3402)(986,3415)(978,3425)
	(969,3433)(959,3440)(949,3446)
	(938,3451)(927,3455)(916,3459)
	(905,3463)(893,3467)(881,3471)
	(869,3475)(856,3479)(842,3484)
	(828,3488)(814,3493)(799,3497)
	(783,3501)(767,3504)(751,3505)
	(735,3504)(713,3498)(693,3487)
	(674,3474)(659,3461)(645,3447)
	(634,3434)(624,3421)(616,3409)
	(609,3398)(602,3386)(596,3374)
	(590,3360)(584,3346)(577,3329)
	(570,3309)(562,3286)(554,3260)
	(547,3230)(541,3197)(537,3162)
	(536,3123)(540,3086)(546,3051)
	(554,3019)(564,2989)(576,2962)
	(589,2936)(603,2912)(618,2889)
	(632,2868)(646,2849)(658,2832)
	(669,2818)(677,2808)(683,2801)
	(686,2797)(688,2795)
\path(728,2752)(729,2750)(732,2746)
	(737,2739)(743,2730)(751,2719)
	(761,2708)(771,2696)(784,2682)
	(799,2668)(816,2653)(837,2637)
	(859,2622)(879,2609)(898,2599)
	(915,2590)(931,2583)(946,2577)
	(959,2572)(971,2568)(979,2565)
	(984,2563)(987,2562)
\path(837,3162)(838,3162)(844,3162)
	(855,3161)(868,3160)(879,3158)
	(888,3156)(895,3153)(900,3149)
	(904,3143)(908,3134)(912,3124)
	(916,3115)(920,3106)(925,3099)
	(931,3095)(942,3091)(956,3089)
	(973,3088)(985,3087)(987,3087)
\path(912,687)(912,686)(912,680)
	(913,669)(914,656)(916,645)
	(918,636)(921,629)(925,624)
	(931,620)(940,616)(950,612)
	(959,608)(968,604)(975,599)
	(979,593)(983,582)(985,568)
	(986,551)(987,539)(987,537)
\path(2037,687)(2037,686)(2037,680)
	(2038,669)(2039,656)(2041,645)
	(2043,636)(2046,629)(2050,624)
	(2056,620)(2065,616)(2075,612)
	(2084,608)(2093,604)(2100,599)
	(2104,593)(2108,582)(2110,568)
	(2111,551)(2112,539)(2112,537)
\path(4287,237)(4287,239)(4286,251)
	(4285,268)(4283,282)(4279,293)
	(4275,299)(4268,305)(4259,309)
	(4249,314)(4240,320)(4231,327)
	(4225,337)(4221,345)(4219,354)
	(4217,367)(4215,382)(4214,401)
	(4213,421)(4213,441)(4212,454)
	(4212,461)(4212,462)
\put(1005.750,1362.000){\arc{1511.673}{1.6952}{2.7334}}
\end{picture}
}
\caption{Irreducible components and other interesting loci in $\cmbar_1(\proj^2,3)$}
\label{cubics}
\end{figure}

The contracted genus $1$ curve could meet one other curve, necessarily
genus $0$ and mapping with degree $3$ (see the upper-left panel of
Fig.~\ref{cubics}).  The general such genus $0$ map will have as image
a nodal cubic.  This component of the moduli space has dimension $10$:
there is an $8$-dimensional family of nodal cubics, plus a
$1$-dimensional choice of where to ``glue'' the elliptic curve, plus a
$1$-dimensional choice of $j$-invariant.  Thus this locus cannot
lie in the closure of the $9$-dimensional main component.

Another possibility is that the contracted genus $1$ curve could meet
{\em two} other curves, one mapping with degree $2$ and one mapping
with degree $1$ (see the upper-right panel of Fig.~\ref{cubics}).  This forms a $9$-dimensional family: a
$2$-dimensional choice for the $2$-pointed elliptic curve ($\dim
\cmbar_{1,2}=2$), plus a $2$-dimensional choice for the image of the
contracted curve in the plane, plus a $4$-dimensional choice of conic
through that point, plus a $1$-dimensional choice of a line through
that point.  Again, for dimensional reasons, all such maps can't lie
in the $9$-dimensional main component.

The final possibility involving a contracted elliptic component is if
the contracted curve meets three other curves, each mapping with
degree $1$ (as lines).  (See the lower-middle panel of
Figure~\ref{cubics}.)  This family has dimension $8$ ($3$ dimensions
for the choice of a point in $\cmbar_{1,3}$, plus a $2$-dimensional
choice of the image of this component in the plane, plus a
$3$-dimensional choice of the three lines through that point).  Thus
there is no dimensional obstruction for all such maps to lie in the
(boundary of the) main component, and indeed they do.

One can extend this analysis to see where the components meet.  
The ``one-tail component'' meets the main component
along the locus of maps where the genus $0$ degree $3$ map has a cusp
precisely where it meets the contracted elliptic curve (see the lower-left
panel of Fig.~\ref{cubics}).
The ``two-tail component'' meets the main component along the locus
of maps where the conic and the line are tangent (see the lower-right panel
of Fig.~\ref{cubics}).
More generally, one can explicitly describe which genus one stable
maps are ``smoothable'' (i.e.\ lie in the main component):

\noindent
{\em Proposition.} 
A genus $1$ stable map $\pi: C \rightarrow \proj^n$ is smoothable
if and only if it is one of two forms:
\begin{itemize}
\item[(i)] $\pi$ contracts no genus $1$ curve, or
\item[(ii)] if $E$ is the maximal connected genus-one curve
contracted by $\pi$, 
and $E$ meets the rest of $C$ (i.e.\ $C' = \overline{C-E}$) at
the points $p_1$, \dots, $p_m$, then $\pi( T_{C', p_1})$, \dots,
$\pi(T_{C, p_m})$ must be a dependent set of vectors in $\pi(T_{\proj^n, \pi(E)})$.
\end{itemize}

This follows readily from the same proof as  \cite[Lemma~5.9]{crelle}.
(More generally, one of the implications holds if $\proj^n$ is replaced by
a smooth target: if $C \rightarrow X$ is smoothable, then one of these
two hold.)

Notice that this proposition ``explains'' the bottom row of
Figure~\ref{cubics}: if $E$ has ``one tail'' ($m=1$), then $\pi(C')$
must have a cusp at that point for the map to be smoothable.  If $E$
has ``two tails'' ($m=2$), then the two branches of $\pi(C')$ must be
tangent at that point for the map to be smoothable.  If $E$ has
``three tails'', then the three branches of $\pi(C')$ must be coplanar
for the map to be smoothable --- but this is automatic in $\proj^2$.

\noindent
{\bf The desingularization.}

We finally describe the desingularization.  We assume $d>0$, as if
$d=0$, $\cmbar_{1,k}(\proj^n,d) \equiv \cmbar_{1,k} \times \proj^n$,
which is already smooth.

Define the {\em $m$-tail locus} of $\cmbar_{1,k}(\proj^n,d)$ to be the
locus maps where there a contracted elliptic curve meets the rest of
the curve and the set of marked points in a total of precisely $m$
points.  (For example, the contracted elliptic curve could contain no
marked points, and meet the rest of the curve in two points; or it
could contain one marked point, and meet the rest of the curve in one
point.)

The $m$-tail locus is the union of a number of components,
which we now describe.  For each $m' \in \{ 1, \dots, d \}$, each
partition $\mu$ of $d$ into $m'$ parts, and each subset $S$ of $\{1,
\dots, k \}$ of size $m-m'$, we have a smooth subvariety (substack,
really) corresponding to maps with a contracted elliptic curve
containing the marked points $S$, and meeting genus $0$ curves mapping
with degrees corresponding to the partition $\mu$.  These may be 
components of $\cmbar_{1,k}(\proj^n,d)$, but may not be (as we saw
in the example of the cubics).

Then the desingularization may be describe as follows:  {\em blow up
the one-tail locus, then the proper transform of the two-tail locus,
etc.}  At each stage, we are blowing up along a smooth center.

We need to blow up in this particular order for the following reason.
Figure~\ref{threefour} shows a map contained in the two-tail,
three-tail, and four-tail locus.  In fact, it is in ``two branches''
of the three-tail locus in the moduli space, corresponding to the two
ways we can select three nodes separating a genus one contracted curve
from the rest of the curve.  Thus the three-tail locus is not smooth
at this point.  Blowing up the two-tail locus will separate these two
branches of the three-tail locus, and the proper transform of the three-tail locus is then smooth (at the points corresponding to this map).

\begin{figure}
\setlength{\unitlength}{0.00083333in}
\begingroup\makeatletter\ifx\SetFigFont\undefined%
\gdef\SetFigFont#1#2#3#4#5{%
  \reset@font\fontsize{#1}{#2pt}%
  \fontfamily{#3}\fontseries{#4}\fontshape{#5}%
  \selectfont}%
\fi\endgroup%
{\renewcommand{\dashlinestretch}{30}
\begin{picture}(4224,1239)(0,-10)
\put(1887,912){\makebox(0,0)[lb]{{\SetFigFont{5}{6.0}{\rmdefault}{\mddefault}{\updefault}genus $1$}}}
\path(1212,912)(1812,312)
\path(3012,912)(2412,312)
\dashline{60.000}(2412,612)(4212,612)
\dashline{60.000}(2412,462)(4212,462)
\dashline{60.000}(1812,612)(12,612)
\dashline{60.000}(12,462)(1812,462)
\path(1362,1212)(1362,762)
\path(1332.000,882.000)(1362.000,762.000)(1392.000,882.000)
\path(2562,12)(2562,462)
\path(2592.000,342.000)(2562.000,462.000)(2532.000,342.000)
\path(2712,1212)(2712,612)
\path(2682.000,732.000)(2712.000,612.000)(2742.000,732.000)
\path(3012,762)(1212,762)
\end{picture}
}
\caption{A map in two branches of the three-tail locus;
the genus one curve is indicated, 
the solid components are the contracted ones, and the arrows represent
a choice of three nodes separating a genus one contracted curve from the rest of the curve}
\label{threefour}
\end{figure}
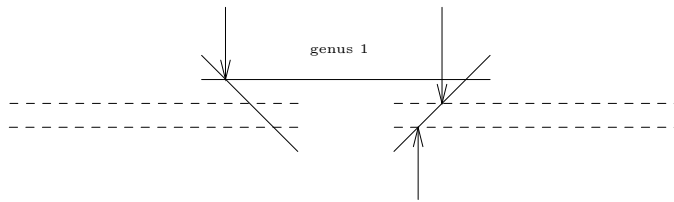

We make a few observations.  

First, this suggests that we should think
of the one-tail locus as the ``most singular locus'', the two-tail
locus as the ``next-most singular locus'', and so on.    This is perhaps
opposite to the order one would expect.

Second, note that blowing up a space (such as
$\cmbar_{1,k}(\proj^n,d)$) may be interpreted as removing the
component (``blowing it out of existence''), and blowing up that
component's scheme-theoretic intersection with the remainder of the
space.  More formally, if $X \cup Y$ is a scheme, with closed subschemes $X$ and $Y$, $\Bl_X (X \cup Y)$ is canonically isomorphic to
$\Bl_{X \cap Y} Y$ by the universal property of blowing up.  Hence we
could equally well describe this construction as blowing up
$\cmbar_{1,k}^0(\proj^n,d)$ along the ``one-tail locus'' of this
space, then the ``two-tail locus'', etc.  (In this case, the first
blow-up, along the one-tail locus, does nothing, as this is already
a Cartier divisor.)  With this interpretation, at each stage we
are still blowing up a space along a smooth center.

For example in the example of cubics, we remove the two non-main
components (the upper-left and upper-right panels of
Figure~\ref{cubics}), blow up the locus corresponding to maps
corresponding to the panel in the lower-right of Figure~\ref{cubics}
(which is a Weil divisor, but not Cartier), then blow up (the proper
transform of) the locus corresponding to the panel in the lower-middle
of Figure~\ref{cubics}.

Third, this construction involves only the underlying curve and the
information of which components are contracted.  By making this precise,
we are led to a candidate definition for more general target
spaces.  Let $\fM_{1,k}$ be the moduli space (Artin stack) of
projective connected, nodal, genus $1$, $k$-pointed nodal curves (over
$\C$).  Construct $\beta: \fM'_{1,k} \rightarrow \fM_{1,k}$ where 
points of $\fM'_{1,k}$ are defined as projective connected nodal genus one
curves with the additional information of a connected union of components of
arithmetic genus $1$ (possibly empty) that is declared to be
contracted.  Then  $\be$ is locally (on the
source) an isomorphism, but is  not separated.  The forgetful morphism
$\cmbar_{1,k}(\proj^n,d) \rightarrow \fM_{1,k}$ naturally factors
through $\fM'_{1,k}$.  If $\widetilde{\fM'}_{1,k}$ is the blow-up of
$\fM'_{1,k}$ along the one-tail locus, then the proper transform of
the two-tail locus, etc., then
$$\cmbar_{1,k}(\proj^n,d) \times_{\fM'_{1,k}}{\widetilde{\fM'}_{1,k}}
$$
 contains $\widetilde{\cm}_{1,k}(\proj^n,d)$
as an irreducible component.  
If $X$ is a complex projective manifold, one can similarly define 
$\widetilde{\cm}_{1,k}(X, \be)$  as the union of components
of
\begin{equation}
\label{bigspace2}
\cmbar_{1,k}(X,\be) \times_{\fM'_{1,k}}{\widetilde{\fM'}_{1,k}}
\end{equation}
generically mapping to $\fM_{1,k}$ (i.e.\ corresponding to maps
with smooth source).
(We have no reasonable modular interpretation
of $\widetilde{\cm}_{1,k}(X, \be)$ in general; taking the
closure is an awkward construction moduli-theoretically.)
Via the exact sequence
for the tangent-obstruction theory of $\cmbar_{1,k}(X,\be)$
in terms of that of $\fM_{1,k}$ and $H^i(C, \pi^* T_X)$, one
can endow \eqref{bigspace2} with a natural virtual fundamental
class.  We expect
this to lead to an algebraic theory of ``reduced genus $1$ Gromov-Witten
invariants'', cf.\ \cite{reduced}.

} 

 \end{document}